\documentclass[12pt]{article}
\begin{document}
\newtheorem{theorem}{Theorem}
\newtheorem{corollary}{Corollary}
\newtheorem{note}{Note}

\begin{center}
\vskip 1cm{\LARGE\bf
Some Integer Sequences \\
Based on Derangements}
\vskip 1cm
\large
Milan Janji\'c \\
Department of Mathematics and Informatics\\
University of Banja Luka\\
Republic of Srpska, Bosnia and Herzegovina\\
mailto:agnus@blic.net
\end{center}

\begin{abstract}
Sequences whose terms are equal to the number of functions with specified properties are
considered. Properties are based on the notion of derangements in a more general sense.
Several sequences which generalize the standard notion of derangements  are thus  obtained.
These sequences generate a number of integer sequences from the well-known Sloane's
encyclopedia.
\end{abstract}
Let $A$ be   an $m\times n$  rectangular area whose elements are from a set $\Omega,$
and let $c_1,\ldots,c_m$ be from $\Omega.$ Following the paper  \cite{mj},  we call
each column of $A$ which is equal to  $[c_1,\ldots,c_m]^T$ an i-column of $A.$
As usual by $[n]$ will be denoted the set $\{1,2,\ldots,n\},$ and by $|X|$ the number of
elements of a finite set $X.$ Mutually disjoint subsets are called blocks. A block with $k$
 elements is called $k$-block. We also denote by $n^{(m)}$ the falling factorials, that is,
$n^{(m)}=n(n-1)\cdots (n-m+1).$  Stirling numbers of the second kind will be denoted by
 $S(m,n).$

We start with the following:
\begin{theorem}\label{t1}
Suppose that $X_1,X_2,\ldots,X_k$ are blocks in  $[m]$ and $Y_1,Y_2,\ldots,Y_k$ are subsets
in $[n].$ Label all functions $f:[m]\to[n]$ by $1,2,\ldots, n^m$ arbitrary and form a
$k\times n^m$ matrix $A=(a_{ij})$ such that $a_{ij}=1$ if $f_j(X_i)\subseteq Y_i,$
and $a_{ij}=0$ otherwise. The number $D_1$ of  i-columns of $A$ consisting of $0$'s is equal
\begin{equation}\label{e1}D=\sum_{I\subseteq [k]}(-1)^{|I|}A(I),\end{equation}
where
\begin{equation}\label{e2}A(I)=n^{|[m]\setminus \cup_{i\in I}X_i|}\cdot \prod_{i\in I}
|Y_i|^{|X_i|},\end{equation}
and $I$ runs over all subsets of $[m]$.
\end{theorem}
\noindent\textbf{Proof.} According to Theorem 1.1 in \cite{mj}, the number $D$ is equal to
the right  side of (\ref{e1}) if $A(I)$ is the maximal number of columns $j$ of $A$ such that
 $a_{ij}=1$ for all $i\in I.$  It follows that $A(I)$ is equal to the number of functions
 $f:[m]\to [n]$ such that  $f(X_i)\subseteq Y_i,\;(i\in I).$ This number is clearly  equal to
 the number on the right side of (\ref{e2}).

In a similar way we obtain the following:
\begin{theorem}\label{t2}
Suppose that $X_1,X_2,\ldots,X_k$ are blocks in $[m]$ and $Y_1,Y_2,\ldots,Y_k$ are subsets of
$[n].$ Label all functions $f:[m]\to[n]$ by $1,2,\ldots, n^m$ arbitrary, and form a
$k\times n^m$ matrix $B=(b_{ij})$ such that $b_{ij}=1$ if $f_j(X_i)=Y_i,$ and $a_{ij}=0$ otherwise.
The number $N$  of i-columns of $A$ consisting of $0$'s is equal
$$N=\sum_{I\subseteq [k]}(-1)^{|I|}B(I),$$
where
$$B(I)=n^{|[m]\setminus \cup_{i\in
I}X_i|}\cdot \prod_{i\in I}|Y_i|!S(|X_i|,|Y_i|),$$ and $I$ runs over all subsets of $[m]$.
\end{theorem}
Depending on the number of elements of $X_1,\ldots,X_k;Y_1,\ldots,Y_k$ it is possible to
obtain a number of different sequences. Consider first the simplest  case when each
$X_1,\ldots,X_k;\;Y_1,\ldots,Y_k$ consists of one element. Then
$$A(I)=n^{m-|I|},$$ so that Theorem 1.2 of \cite{mj} may be applied. We thus obtain the following
 consequence of Theorem \ref{t1}.
 \begin{corollary}
 Given distinct $x_1,\ldots,x_k$ in $[m]$ and arbitrary  $y_1,\ldots,y_k$ in $[n]$, then the
 number $D_{11}(m,n,k)$ of functions $f:[m]\to[n]$ such that $$f(x_i)\not=y_i,\;(i=1,2,\ldots,k),$$
is equal
$$D_{11}(m,n,k)=\sum_{i=0}^k(-1)^i{k\choose i}n^{m-i}\left(=n^{m-k}(n-1)^k\right).$$
\end{corollary}

A number of sequences in \cite{ns} is generated by this simple function. Some of them are stated
in the following:
\begin{center}\textbf{Table 1.}\end{center}
\begin{center}
\begin{tabular}{ll}
1.  $A001477(n)=D_{11}(1,n,1)$,&
2.  $A002378(n)=D_{11}(2,n,1)$,\\
3.  $A045991(n)=D_{11}(3,n,1),$&
4.  $A085537(n)=D_{11}(4,n,1)$,\\
5.  $A085538(n)=D_{11}(5,n,1),$&
6.  $A085539(n)=D_{11}(6,n,1),$\\
7.  $A000079(n)=D_{11}(n,2,1),$&
8.  $A008776(n)=D_{11}(n,3,1),$\\
9.  $A002001(n)=D_{11}(n,4,1),$&
10. $A005054(n)=D_{11}(n,5,1),$\\
11. $A052934(n)=D_{11}(n,6,1),$&
12. $A055272(n)=D_{11}(n,7,1),$\\
13. $A055274(n)=D_{11}(n,8,1),$&
14. $A055275(n)=D_{11}(n,9,1),$\\
15. $A052268(n)=D_{11}(n,10,1),$&
16. $A055276(n)=D_{11}(n,11,1),$\\
17. $A000290(n)=D_{11}(2,n,2),$&
18. $A011379(n)=D_{11}(3,n,2),$\\
19. $A035287(n)=D_{11}(4,n,2),$&
20. $A099762(n)=D_{11}(5,n,2),$\\
21. $A000079(n)=D_{11}(n,2,2),$&
22. $A003946(n)=D_{11}(n,3,2),$\\
23. $A002063(n)=D_{11}(n,4,2),$&
24. $A055842(n)=D_{11}(n,5,2),$\\
25. $A055846(n)=D_{11}(n,6,2),$&
26. $A055270(n)=D_{11}(n,7,2),$\\
27. $A055847(n)=D_{11}(n,8,2),$&
28. $A055995(n)=D_{11}(n,9,2),$\\
29. $A055996(n)=D_{11}(n,10,2),$&
30. $A056002(n)=D_{11}(n,11,2),$\\
31. $A056116(n)=D_{11}(n,12,2),$&
32. $A076728(n)=D_{11}(n,n,2),$\\
33. $A000578(n)=D_{11}(3,n,3),$&
34. $A005051(n)=D_{11}(n,3,3),$\\
35. $A056120(n)=D_{11}(n,4,3),$&
36. $A000583(n)=D_{11}(4,n,4),$\\
37. $A101362(n)=D_{11}(5,n,4),$&
38. $A118265(n)=D_{11}(n,4,4).$
\end{tabular}
\end{center}

 Suppose that $$|X_1|=|X_2|=\ldots=|X_k|=1,\;|Y_1|=|Y_2|=\cdots=|Y_k|=2.$$ Then
$$A(I)=2^in^{m-|I|}.$$  We  may again  apply Theorem 1.2 in [2] to obtain
the following:
\begin{corollary}
Given distinct $x_1,\ldots,x_k$ in $[m]$ and arbitrary 2-sets  $Y_1,\ldots,Y_k$
 in $[n]$, then the number $D_{12}(m,n,k)$ of functions $f:[m]\to[n]$ such that
 $$f(x_i)\not\in Y_i,\;(i=1,2,\ldots,k),$$
is equal
$$D_{12}(m,n,k)=\sum_{i=0}^k(-2)^i{k\choose i}n^{m-i}\left(=n^{m-k}(n-2)^k\right).$$
\end{corollary}

This function also generates a number of sequences in \cite{ns}.The following table
contains some of them.
\begin{center}\textbf{Table 2.}\end{center}
\begin{center}
\begin{tabular}{ll}
1. $A000027(n)=D_{12}(1,n,1),$&
2. $A005563(n)=D_{12}(2,n,1)$\\
3. $A027620(n)=D_{12}(3,n,1),$&
4. $A000244(n)=D_{12}(n,3,1)$,\\
5. $A004171(n)=D_{12}(n,4,1),$&
6. $A005053(n)=D_{12}(n,5,1),$\\
7. $A067411(n)=D_{12}(n,6,1),$&
8. $A000290(n)=D_{12}(2,n,2),$\\
9. $A0002444(n)=D_{12}(n,3,2),$&
10. $ A000578(n)=D_{12}(3,n,3),$\\
 11. $ A081294(n)=D_{12}(n,4,3),$&
 12. $A000583(n)=D_{12}(4,n,4),$
\end{tabular}
\end{center}

If, in the conditions of Theorem\ref{t1}, hold
 $$|X_1|=\cdots=|X_k|=2;\;|Y_1|=\cdots|Y_k|=1,$$
then $$A(I)=n^{m-2|I|},$$ so that we have the following:
\begin{corollary}
Suppose that  $X_1,\ldots,X_k$ are 2-blocks in $[m],$  and  $y_1,\ldots,y_k$
 arbitrary elements  in $[n]$, then the number $D_{21}(m,n,k)$ of functions $f:[m]\to[n]$
 such that
 $$f(X_i)\not=\{y_i\},\;(i=1,2,\ldots,k)$$
is equal
$$D_{21}(m,n,k)=\sum_{i=0}^k(-1)^i{k\choose i}n^{m-2i}\left(=n^{m-2k}(n^2-1)^k\right).$$
\end{corollary}

We also state some sequences in \cite{ns} generated by this function.
\begin{center}\textbf{Table 3.}\end{center}
\begin{center}
\begin{tabular}{ll}
1. $A005563(n)=D_{21}(2,n,1),$&
2. $A007531(n)=D_{21}(3,n,1)$\\
3. $A047982(n)=D_{21}(4,n,1),$&
4. $A005051(n)=D_{21}(n,3,1)$,\\
5. $A005010(n)=D_{21}(n,2,2),$&
\end{tabular}
\end{center}
Take finally the case  $|X_i|=|Y_i|=2,\;(i=1,2,\ldots,k)).$ We have now
$$A(I)=4^{|I|}\cdot n^{m-2|I|}.$$ We thus obtain the following consequence of Theorem\ref{t1}.
\begin{corollary}
 Let  $X_1,\ldots,X_k$ in $[m]$ be 2-blocks,  and  $Y_1,\ldots,Y_k$
 in $[n]$ be  arbitrary 2-sets. Then the number $D_{22}(m,n,k)$ of functions $f:[m]\to[n]$
 such that
 $$f(X_i)\not\subset Y_i,\;(i=1,2,\ldots,k)$$
is equal
$$D_{22}(m,n,k)=\sum_{i=0}^k(-4)^i{k\choose i}n^{m-2i}\left(=n^{m-2k}(n^2-4)^k\right).$$
\end{corollary}

A few sequences in \cite{ns}, given in the next table, is defined by this function.

\begin{center}\textbf{Table 4.}\end{center}
\begin{center}
\begin{tabular}{ll}
1. $A005030(n)=D_{22}(n,3,1),$&
2. $A002001(n)=D_{22}(n,4,1)$\\
3. $A002063(n)=D_{22}(n,4,2),$&\\
\end{tabular}
\end{center}

Take now the case  $|X_i|=|Y_i|=2,\;(i=1,2,\ldots,k))$ in the conditions of Theorem \ref{t2}.
We have $$B(I)=2^{|I|}\cdot n^{m-2|I|}.$$ Thus we have the next:
\begin{corollary}
 Let  $X_1,\ldots,X_k$ be 2-blocks in $[m]$ and  $Y_1,\ldots,Y_k$ arbitrary 2-sets
 in $[n].$ Then the number $S_{22}(m,n,k)$ of functions $f:[m]\to[n]$
 such that
 $$f(X_i)\not= Y_i,\;(i=1,2,\ldots,k)$$
is equal
$$S_{22}(m,n,k)=\sum_{i=0}^k(-2)^i{k\choose i}n^{m-2i}\left(=n^{m-2k}(n^2-2)^k\right).$$
\end{corollary}

The sequence A005032 in \cite{ns} is generated by this function.

We shall now consider injective functions from $[m]$ to $[n],\;(m\leq n).$
We start with the following:
\begin{theorem}
Let $X_1,X_2,\ldots,X_k$ be  blocks in $[m]$ and $Y_1,Y_2,\ldots,Y_k$  blocks in $[n]$
 such that
$$|X_i|=|Y_i|,\;(i=1,2,\ldots,k).$$
 If a $k\times n^{(m)}$ matrix $A$ is defined such that $a_{ij}=1$ if $f_j(X_i)=Y_i$ and
 $a_{ij}=0$ otherwise, then the number $I(m,n,k)$ of i-columns of $A$ consisting of $0$'s
is equal
$$I_k(m,n)=\sum_{I\subseteq [k]}(-1)^{|I|}(n-|\cup_{i\in I}X_i|)^{(m-|\cup_{i\in
I}X_i|)}\cdot \prod_{i\in I}|X_i|!.$$
\end{theorem}
\textbf{Proof.} In this case we have
$$A(I)=(n-|\cup_{i\in I}X_i|)^{(m-|\cup_{i\in I}X_i|)}\cdot\prod_{i\in I}|X_i|!,$$
so that theorem follows from Theorem 1.1. in \cite{mj}.

We shall also state some particular cases of this theorem. Suppose first that
$$|X_i|=|Y_i|=1,\;(i=1,\ldots,k).$$
The number $A(I)$ in this case is equal
$$(n-|I|)^{(m-|I|)}.$$  We thus obtain the following:
\begin{corollary}
For disjoint $x_1,\ldots,x_k$ in $[m]$ and disjoint $y_1,\ldots,y_k$ in $[m],$  the number
$I_1(m,n,k)$ of injections $f:[m]\to [n]$ such that
$$f(x_i)\not=y_i,\;(i=1,2,\ldots,k)$$
is equal
$$I_1(m,n,k)=\sum_{i=0}^k(-1)^i{k\choose i}(n-i)^{(m-i)}.$$
\end{corollary}

\begin{note}

Since obviously holds $D(n)=I(n,n,n),$ where $D(n)$ is the number of derangements
of $n$ elements, this function is an extension of derangements.
\end{note}

 There are a number of sequences in \cite{ns} that are generated by this function.
 We state  some of them in the next table.
\begin{center}\textbf{Table 5.}\end{center}
\begin{center}
\begin{tabular}{ll}
1. $A000290(n)=I(2,n,1),$&
2. $A045991(n)=I(3,n,1)$\\
3. $A114436(n)=I(3,n,1)$&
4. $A047929(n)=I(4,n,1),$\\
5. $A001563(n)=I(n,n,1),$&
6. $A001564(n)=I(n,n,2),$\\
7. $A001565(n)=I(n,n,3),$&
8. $A002061(n)=I(2,n,2)$,\\
9. $A027444(n)=I(3,n,2),$&
10. $A058895(n)=I(4,n,2)$,\\
11. $A027444(n)=I(3,n,2),$&
12. $A074143(n)=I(n-1,n,1)$,\\
13. $A001563(n)=I(n-1,n,1),$&
14. $A094304(n)=I(n-1,n,1)$,\\
15. $A109074(n)=I(n-1,n,1),$&
16. $A094258(n)=I(n-1,n,1)$,\\
17. $A001564(n)=I(n-1,n,2),$&
18. $A001565(n)=I(n-1,n,3)$,\\
19. $A001688(n)=I(n-1,n,4)$&
20. $A001689(n)=I(n-1,n,5)$,\\
21. $A023043(n)=I(n-1,n,6),$&
22. $A023044(n)=I(n-1,n,7)$,\\
23. $A023045(n)=I(n-1,n,8),$&
24. $A023046(n)=I(n-1,n,9),$\\
25. $A023407(n)=I(n-1,n,10),$&
26. $A001563(n)=I(n-2,n,1),$\\
27. $A001564(n)=I(n-2,n,2),$&
28. $A061079(n)=I(n,2n,1).$
\end{tabular}
\end{center}

As a  special case we also have  the following generalization of derangements.

\begin{corollary}
If $X_1,X_2,\ldots,X_n$ is a partition of $[kn]$
such that $$|X_i|=k,\;(i=1,2,\ldots,n)),$$
then the number $D(n,k)$ of permutations $f$ of $[kn]$ such that $f(X_i)\not=X_i,\;(i=1,2,\ldots,n)$ is equal
$$D(n,k)=\sum_{i=0}^n(-1)^i(k!)^i(nk-ik)!.$$
\end{corollary}

For $k=1$ we obtain the standard formula for derangements.
\begin{note}
From the preceding formula the following sequences in \cite{ns} are derived:\\
$A128805,  A127888, A116221, A116220, A116219.$
\end{note}

\bigskip
\hrule
\bigskip

\noindent 2000 {\it Mathematics Subject Classification}:
Primary 05A10.

\noindent \emph{Keywords: }
derangements, Stirling numbers of the second kind.

\bigskip
\hrule
\bigskip

\noindent Concerned with sequences:\\
A001477, A002378, A045991, A085537, A085538, A085539, A000079,\\
A008776, A002001, A005054, A052934, A055272, A055274, A055275,\\
A052268, A055276, A000290, A011379, A035287, A099762, A000079,\\
A003946, A002063, A055842, A055846, A055270, A055847, A055995,\\
A055996, A056002, A056116, A076728, A000578, A005051, A056120,\\
A000583, A101362, A118265, A000027, A005563, A027620, A000244,\\
A004171, A005053, A067411, A000290, A002444, A000578, A081294,\\
A000583, A005563, A007531, A047982, A005051, A005010, A005032,\\
A005030, A002001, A002063, A005032, A000290, A045991, A114436,\\
A047929, A001563, A001564, A001565, A002061, A027444, A058895,\\
A027444, A074143, A001563, A094304, A109074, A094258, A001564,\\
A001565, A001688, A001689, A023043, A023044, A023045, A023046,\\
A023407, A001563, A001564, A061079, A128805,  A127888, A116221,\\
A116220, A116219.

\end{document}